\documentclass[11pt,a4paper]{article}
\usepackage[utf8]{inputenc}
\usepackage{amsmath}
\usepackage{amsfonts}
\usepackage{amssymb}
\usepackage{graphicx}
\usepackage[left=3.5cm,right=2.5cm,top=2cm,bottom=2cm]{geometry}
\usepackage{authblk}

\title{\textbf{Common Fixed Point Theorems in Fuzzy Metric Space with Applications }}
\author[]{Rachana Soni\thanks{rachanasoni007@gmail.com}}
\affil[]{Department of Applied Mathematics,\\ Bhilai Institute of Technology, Durg (C.G)India.}
 \begin{document}
\date{}
\maketitle
\begin{abstract}\textbf{ In this paper, we introduce a new class of implicit function to prove common fixed point theorems in fuzzy metric space. Moreover we define a new altering distance in terms of integral and utilize the same to deduce integral type contractive conditions. Secondly we present  application of main results  to the system of functional equations. At the end we give an example in support of results of the paper.}
\end{abstract}
\textbf{Subject:} 2010 AMS Classification: 47H10, 54H25.\\ \vspace{2mm}

\section{Introduction} \hspace {5mm} The study of fuzzy mathematics began to explore after Zadeh \cite{Z} introduced the idea of fuzzy sets in 1965 to encounter ambiguity of our day to day life. Many authors have worked out on the theory of fuzzy sets and its applications and explored it successfully. Nevertheless when the vagueness is due to fuzziness instead of randomness, as in the measurement of distance, the concept of fuzzy metric space sounds more appropriate. There are numerous definitions of fuzzy metric space (\cite{G2,G1,KM}). George and Veeramani modified the concept of fuzzy metric space initiated by Kramosil and Michalek. In \cite{G2} they explained that every metric genrates a fuzzy metric. This has recently found very significant applications in quantum particle physics particularly in connection with both string and $\epsilon^{\infty}$ theory (see \cite{MS}). \\

It is well known that fixed point theory is one of the most researched areas in Nonlinear Analysis. It can been applied to very different abstract metric spaces and, in particular, recently, many fixed point results have been established in the setting of fuzzy metric spaces \\

Now, if we talk about implicit functions, Popa \cite{VP2,VP1} introduced the idea of implicit function to prove a common fixed point theorem in metric spaces. Imdad and Ali \cite{IA} generalizes the result of Popa in fuzzy metric spaces. Consequently implicit relations are used as a tool for finding common fixed point of contraction maps. Recently J. Ali et al. \cite{JAli} set up great piece of work in such manner via $\delta$-distances in the settings of semi-metric spaces.\\

In this paper, we further attempt to establish common fixed point theorems involving implicit function and altering distances. As first application we produce integral type contractive conditions
from new class of altering distances in the settings of fuzzy metric space which is being introduced by us. As second application of main results we present solution of system of functional equations by using main results. Lastly we provide an example to validate our main results.

\section{Mathematical Preliminaries}
 \textbf{Definition 2.1} (\cite{Z}) Let X be any set. A fuzzy set A in X is a function with domain X and values in [0, 1]. \\

\noindent \textbf{Definition 2.2.}(\cite{KM})  A binary operation $\ast : [0,1] \times [0,1] \rightarrow [0,1] $ is continuous t-norm if $\ast$ satisfies the following conditions: \\ 
(i) \ $\ast$ is commutative and associative, \\
(ii) $\ast$ is continuous, \\
(iii)\ a $\ast$1= a for every a $\in$ [0,1],  \\
(iv)\ a $\ast$ b $\leq$ c $\ast$ d whenever a $\leq$ b and c $\leq$ d for all a, b, c, d $\in$ [0,1]. \\

\noindent\textbf{Definition 2.3.}(\cite{G2}) The 3-tuple (X,M,$\ast$) is called a fuzzy metric space(FM-space) if X is an arbitrary set $\ast$ is a continuous t-norm and M is a fuzzy set in $X^{2}	\times (0,\infty)$ satisfying, for every $x, y, z \in X$ and $t, s > 0,$ the following conditions: \\
(FM-1)$ M(x, y, 0) = 0,$\\ 
(FM-2) $M(x, y, t) = 1$ for all $t > 0$ if and only if $x = y$,\\
(FM-3) $M(x, y, t) = M(y, x, t),$\\
(FM-4)$ M(x, y, t) \ast M(y, z, s) \leq M(x, z, t+s),$\\
(FM-5)$ M(x, y , \cdot): (0, \infty) \rightarrow [0, 1]$ is continuous.\\

\noindent \textbf{Remark 1} (\cite{GR}). Let $(X, M, \ast)$ be a fuzzy metric space. Then $M(x, y, \cdot)$ is nondecreasing on $(0, \infty)$ for all $x,y \in X.$\\

\noindent \textbf{Remark 2} (\cite{G2}) Let $(X, M, \ast)$ be a fuzzy metric space. Then $M(x, y, \cdot)$ is continuous function on $X^2 \times (0, \infty)$ for all $x,y \in X$.\\

\noindent \textbf{Definition 2.4} (\cite{JG}) A pair of self mappings (A,B) of a fuzzy metric space (X,M,$\ast$) is said to be commuting if M(ABx, BAx, t) =1 for all x $\in$ X.\\

 Pant \cite{P} has given thought of R-weakly commuting maps in metric spaces in 1994. Vasuki \cite{RV} practiced the term R-weakly commuting mapping in fuzzy metric spaces and proved some common fixed point theorems for these mappings.\\

\noindent\textbf{Definition 2.5.} (\cite{RV}) The pair of self maps (A, B) of a fuzzy metric space (X, M , $\ast$) is said to be, \\
(1) weakly commuting if $M(ABx, BAx, t) \geq M(Ax, Bx, t)$ for all $x \in X \ and \ t > 0,$\\
(2) R-weakly if there exists $R > 0$ such that $M(ABx, BAx, t) \geq M(Ax, Bx, t/R)$ for all $x \in X$ and $t \geq 0$. \\

\noindent\textbf{Remark 3} (\cite{SNM}) Let $(X, M, \ast )$ be a fuzzy metric space. If there exists $r \in (0,1)$ such that $M(x, y, rt) \geq M(x, y, t)$ for all $x,y \in X$ and $t > 0$, then $x = y$\\

\noindent\textbf{Definition 2.6.} (\cite{J})  A pair of self mappings (A, B) of a fuzzy metric space $(X,M,\ast)$ is said to be compatible (or asymptotically commuting) if for all $t > 0$\\

$\lim_{n\to \infty} M(ABx_{n}, BAx_{n}, t) = 1, $\\

 whenever $\lbrace x_{n} \rbrace$ is a sequence in X such that  $\lim_{n\to \infty} Ax_{n} = \lim_{n\to \infty} Bx_{n} =u$ for some u $\in X. $ Also the pair (A, B) is called noncompatible, if there exists a sequence ${x_n}$ in X such that $\lim_{n\to \infty}Ax_{n} = \lim_{n\to \infty} Bx_{n} =u$, but either $\lim_{n\to \infty} M (ABx_n, BAx_n, t) \neq 1 $ or the limit does not exist.  \\
\noindent\textbf{Definition 2.7.} (\cite{PA}) A pair of self mappings $(A, B)$ of a fuzzy metric space $(X, M, \ast)$ is said to be weakly compatible if they commute at the coincidence points i.e., if $Au=Bu$ for some $u= X,$ then $ABu=BAu.$ \\
\noindent\textbf{Definition 2.8.} (\cite{MI}) A pair of self mappings $(A, B)$ of a fuzzy metric space $(X, M, \ast)$ is said to have the property (E.A.) if there exists a sequence $\lbrace x_n \rbrace $ in X such that $lim_{n\longrightarrow \infty} Ax_n = lim_{n\longrightarrow \infty} Bx_n = y $ for some $y \in X.$ \\

\noindent We can see that compatible as well as noncompatible pairs satisfy the property (E.A.). \\

\noindent\textbf{Definition 2.9.} (\cite{MI}) Two pairs of self mappings $(A, P)$ and $(B, Q)$ defined on fuzzy metric space $(X, M, \ast$) are said to share common property (E.A.) if there exist sequences $\lbrace x_n\rbrace$ and $\lbrace y_n \rbrace$ in X such that,\\ $\lim_{n\to \infty} Ax_n = lim_{n\to \infty} Sx_n = lim_{n\to \infty} By_n = lim_{n\to \infty} Qy_n = z$ for some $z \in X.$ \\

Pathak et al. \cite{HKP} improved the notion of R-weakly commuting mappings in metric spaces by introducing the notions of R-weakly commutativity of
type $(A_g)$ and R-weakly commutativity of type $(A_f )$. Imdad and Ali 
enhanced the notion of R-weakly commutativity of type $(A_g)$ and R-weakly com-
mutativity of type $(A_f )$ in fuzzy metric space
and after that they  launched the concept of R-weakly commuting maps of type
$(P)$ in fuzzy metric spaces in 2008.\\

\noindent \textbf{Definition 2.10} A pair of self-mappings (A, S) of a fuzzy metric space $(X,M, \ast)$ is
said to be\\

(i) R-weakly commuting mappings of type $(A_g)$ if there exists some $R > 0$ such
that $M(SAx, AAx, t) \geq M(Ax, Sx, \frac{t}{R})$ 
 for all $x \in X$ and $t > 0$.\\
 
(ii) R-weakly commuting mappings of type $(A_f )$ if there exists some $R>0$ such that
$M(ASx, SSx, t) \geq M(Ax, Sx, \frac{t}{R})$, for all $x \in X$ and $t > 0.$\\

(iii) R-weakly commuting mappings of type (P) if there exists some $R > 0$ such
that $M(AAx, SSx, t) \geq M(Ax, Sx, \frac{t}{R})$, for all $x \in X$ and $t > 0.$

\section{Implicit Functions} Popa gave the concept of implicit functions in fixed point theory. Various authors utilized this idea to prove fixed point theorems. To define our implicit function, let $\Psi$ be the family of all continuous functions $\mathcal{\psi} : [0, 1]^4\longrightarrow \mathfrak{R} $ satisfying following conditions.\\
$(\mathcal{\psi}_1):\mathcal{\psi}$  is nondecreasing in first argument, \\
$(\mathcal{\psi}_2): \mathcal{\psi}(u,0,u,0) \geq 0 \ \Rightarrow \ u \geq 0$, \\
$(\mathcal{\psi}_3): \mathcal{\psi}(u,0,0,u) \geq 0 \ \Rightarrow \ u \ \geq 0$, \\
$(\mathcal{\psi}_4): \mathcal{\psi}(u,u,0,0) \geq 0 \ \Rightarrow \ u \geq 0$. \\

\noindent \textbf{Example 2.1.} Define $\mathcal{\psi}(u_1,u_2,u_3,u_4): [0,1]^4 \longrightarrow \mathfrak{R}$ as \\

$\psi(u_1,u_2,u_3,u_4) = u_1 - \delta(max\lbrace u_2,u_3,u_4 \rbrace)$\\ 

where $\delta : [0,1] \longrightarrow [0,1]$ is a upper semicontinuous function such that $\delta(u) < u$ for $u > 0.$\\

\noindent \textbf{Example 2.2.} Define $\psi(u_1,u_2,u_3,u_4) : [0,1]^4 \longrightarrow \mathfrak{R}$ as \\

$\psi(u_1,u_2,u_3,u_4) = u_1 - k \ min\lbrace u_2,u_3,u_4 \rbrace$ where $0 < k < 1$. \\

\noindent \textbf{Example 2.3.} Define $\psi(u_1,u_2,u_3,u_4) : [0,1]^4 \longrightarrow \mathfrak{R}$ as \\

$\psi(u_1,u_2,u_3,u_4) =  u_1 - \delta(u_2,u_3,u_4) $\\

\noindent where $\delta: [0,1]^3 \to [0,1]$ is upper semicontinuous function such that max$\lbrace \delta( 0,u,0),\delta(0,0,u),\delta(u,0,0)\rbrace <u  $ for every $u>0.$\\

\noindent \textbf{Example 2.4.} Define $\psi(u_1,u_2,u_3,u_4) : [0,1]^4 \longrightarrow \mathfrak{R}$ as \\

$\psi(u_1,u_2,u_3,u_4) = u_1 -ku_2 - min \lbrace u_3, u_4 \rbrace$ where $0 < k < 1$. \\

Let  $\Phi = \lbrace \phi$  is a Lebesgue integrable mapping such that  $\int_0^\epsilon \phi(u)du > 0 \ for \ all \ \epsilon > 0 \rbrace$. Now, we give examples  of integral type functions\\

\noindent \textbf{Example 2.5.} Define $\psi(u_1,u_2,u_3,u_4) : [0,1]^4 \longrightarrow \mathfrak{R}$ as \\

$$\psi(u_1,u_2,u_3,u_4) = \int_0^{1-u_1} \phi(x)dx - a \ max \Big\lbrace \int_0^{1-u_2} \phi(x)dx, \int_0^{1-u_3} \phi(x)dx, \int_0^{1-u_4}\phi(x)dx \Big\rbrace $$ \\
where   $0 \leq a < 1$ and $\phi \in \Phi$. \\

\noindent \textbf{Example 2.6.} Define $\psi(u_1,u_2,u_3,u_4) : [0,1]^4 \longrightarrow \mathfrak{R}$ as \\

$$\psi(u_1,u_2,u_3,u_4) = \int_0^{1-u_1} \phi(x) dx - \delta \Big(max\Big\lbrace \int_0^{1-u_2}\phi(x)dx ,\int_0^{1-u_3}\phi(x)dx, \int_0^{1-u_4}\phi(x)dx\Big\rbrace \Big)$$\\

\noindent where $\delta:  \mathfrak{R}_+ \longrightarrow \mathfrak{R}_+ $ is continuous function such that $\delta(u) < u$ for all $u > 0$ and $\phi \in \Phi$. \\

Verification of conditions $(\psi_1), (\psi_2), (\psi_3)$ and $(\psi_4)$ in setting of Examples 2.1 - 2.6 are obvious.

\section{Main results}

\textbf{Definition 4.1} A function $\varphi: [0,1] \longrightarrow [0,1]$ is called an altering distance function \cite{SY} if it satisfies the followings; \\
(ad1) $\varphi$  is strictly decreasing and continuous, \\
(ad2) $\varphi (\lambda) = 0$ if and only if $\lambda = 1.$ \\

\noindent \textbf{Theorem 4.1} Let $A, B, F \ and \ G$ be self mappings of a fuzzy metric space $(X, M, \ast)$. Suppose that \\

\noindent (a) Either the pair (A, F) or (B, G) share the property (E.A.),\\
(b) Either B(X) $\subset$ F(X) or G(X) $\subset$ A(X),\\
(c) A(X) is a closed subset of X(or B(X) is  closed subset of X),\\
(d) $\mathcal{\psi} \in \Psi$ and $\varphi$ is an altering distance function for all $x,y \ \in \ X,$ such that \\

$\psi \Big( \varphi(M(Fx, Gy, t)), \varphi(M(Ax, By, t)), \varphi(M(Ax, Fx, t)), \varphi(M(By, Gy, t)) \Big) \geq 0. $ (4.1.1) \\

\noindent Then  pairs (A, F) and (B, G) have a common coincidence point. Mappings $A, B, F$ and $G$ have a unique  common fixed point if  pairs (A, f) and (B, G) are weakly compatible. \\

\noindent \emph{Proof }The pair (A, F) share property (E.A.), then there exists a sequence $\lbrace x_n\rbrace$ in X such that \\
$\lim_{n \to \infty} Ax_n = \lim_{n \to \infty} Fx_n = z,$ for some $z \in X.$ \\

\noindent Since F(X) $\subset$ B(X), there exists $\lbrace y_n \rbrace$ for each $\lbrace x_n \rbrace$ such that $By_n = Fx_n.$ Thus $\lim_{n \to \infty}By_n = \lim_{n \to \infty}Fx_n = z.$ Now we claim that $\lim_{n \to \infty} Gy_n = z$, from (4.1.1) we have \\

 $\psi \Big( \varphi(M(Fx_n, Gy_n, t)), \varphi(M(Ax_n, By_n, t)), \varphi(M(Ax_n, Fx_n, t)), \varphi(M(By_n, Ty_n, t)) \Big) \geq 0$. \\

 $\psi \Big( \varphi(M(Fx_n, Gy_n, t)), \varphi(1), \varphi(1), \varphi(M(Fx_n, Gy_n, t)) \Big) \geq 0$  (as $n \to \infty$) \\

 \noindent From $\mathcal{\psi}_3$ of implicit functions and (ad2) of altering distance functions we get\\

$M(Fx_n, Gy_n, t) \geq 1 $ thus $Fx_n = Gy_n$ implies $\lim_{n \to \infty}Fx_n = \lim_{n \to \infty}Gy_n= z $ \\

Now, we see that \\ 

$\lim_{n \to \infty}Ax_n =  \lim_{n \to \infty}Fx_n = \lim_{n \to \infty}By_n= \lim_{n \to \infty}Gy_n = z $ for $z \in$ X.\\

\noindent If A(X) is a closed subset of X, then $\lim_{n \to \infty}Ax_n = z \in A(X),$ thus there exists a point $v \in X$ such that $Av = z$. We now claim that $Av = Fv$, we put $x = v$ and $y = y_n$  in the condition (4.1.1) \\

 $\psi \Big( \varphi(M(Fv, Gy_n, t)), \varphi(M(Av, By_n, t)), \varphi(M(Av, Fv, t)), \varphi(M(By_n, Gy_n, t)) \Big) \geq 0. $\\

$\psi \Big( \varphi(M(Fv, z, t)), \varphi(M(z, z, t)), \varphi(M(z, Fv, t)), \varphi(M(z, z, t)) \Big) \geq 0. $ (as $n \to \infty$) \\

$\psi \Big( \varphi(M(Fv, z, t)), \varphi(1), \varphi(M(Fv, z, t)), \varphi(1) \Big) \geq 0. $ (as $n \to \infty$)\\

From $(\mathcal{\psi}_2)$ and (ad2) we get\\

$Fv = z = Av$. We get $v$ is point of coincident of self mappings $(A, F).$  \\

Since $F(X) \subset B(X)$ and $Fv \in S(X)$, there will be existence of the point $s \in X$ such that $Bs = Fv = z$. We now claim that $Bs = Gs$, again using (4.1.1), we put $x = v$ and $y = s$  \\

 $\psi \Big( \varphi(M(Fv, Gs, t)), \varphi(M(Av, Bs, t)), \varphi(M(Av, Fv, t)), \varphi(M(Bs, Gs, t)) \Big) \geq 0. $ \\

 $\psi \Big( \varphi(M(Bs, Gs, t)), \varphi(M(z, z, t)), \varphi(M(z, z, t)), \varphi(M(Bs, Gs, t)) \Big) \geq 0. $\\

 $\psi \Big( \varphi(M(Bs, Gs, t)), \varphi(1), \varphi(1), \varphi(M(Bs, Gs, t)) \Big) \geq 0. $\\
 
 From $(\mathcal{\psi}_3)$ and (ad2) we get\\
 
$ Bs = Gs = z$\\

This shows that coincidence point of the self mappings $(B, G)$ is s.\\

 The pair of self mappings $(A, F)$ is weakly compatible, there must exist a point $v$ where the pair of self mappings (A, F) commutes that is $Az = AFv = FAv = Fz$. Now we claim that z is a common fixed point of the pair $(A, F),$ we put $x = z$ and $y = s$ in (4.1.1)\\

 $\psi \Big( \varphi(M(Fz, Gs, t)), \varphi(M(Az, Bs, t)), \varphi(M(Az, Fz, t)), \varphi(M(Bs, Gs, t)) \Big) \geq 0. $ \\ 
 
  $\psi \Big( \varphi(M(Az, z, t)), \varphi(M(Az, z, t)), \varphi(M(Az, Az, t)), \varphi(M(z, z, t)) \Big) \geq 0. $ \\ 
 
  $\psi \Big( \varphi(M(Az, z, t)), \varphi(M(Az, z, t)), \varphi(1), \varphi(1) \Big) \geq 0. $ \\ 
 
 From $(\mathcal{\psi}_4)$ and (ad2) we get \\

 $Az = z = Fz$, which implies that z is the common fixed point of the pair $(A, F)$. \\
 
 The pair $(B, G)$ is also weakly compatible, then there exists a point s at which pair $(B, G)$ commutes that is $Bz = BGs = GBs = Gz.$ Now we put $x = v$ and $y = z$ in (4.1.1) \\

$\psi \Big( \varphi(M(Fv, Gz, t)), \varphi(M(Av, Bz, t)), \varphi(M(Av, Fv, t)), \varphi(M(Bz, Gz, t)) \Big) \geq 0. $ \\

$\psi \Big( \varphi(M(z, Bz, t)), \varphi(M(z, Bz, t)), \varphi(M(z, z, t)), \varphi(M(Bz, Bz, t)) \Big) \geq 0. $ \\

$\psi \Big( \varphi(M(z, Bz, t)), \varphi(M(z, Bz, t)), \varphi(1), \varphi(1) \Big) \geq 0. $\\

From $(\mathcal{\psi}_4)$ and (ad2)we get \\

$Bz = z = Gz$ which shows that z is the common fixed point of the pair $(B, G)$. Thus we proved that z is common fixed point of both the pairs $(A, F)$ and $(B, G)$.\\

\emph{Uniqueness:} Uniqueness is easy outcome of inequality (4.1.1) in vision of $(\mathcal{\psi}_4).$  \\
This completes the proof.\\

\noindent \textbf{Corollary 4.1} Let $A \ and \ F$ be self mappings of a fuzzy metric space $(X, M, \ast)$. Suppose that \\

\noindent (a) The pair $(A, F)$ share the property (E.A.),\\
(b) $A(X)$ is a closed subset of $X$ and \\
(d) $\mathcal{\psi} \in \Psi$ and $\varphi$ is an altering distance function for all $x,y \ \in \ X,$ such that \\

 $\psi \Big( \varphi(M(Fx, Fy, t)), \varphi(M(Ax, Ay, t)), \varphi(M(Ax, Fx, t)), \varphi(M(Ay, Fy, t)) \Big) \geq 0. $  \\

\noindent Then the pair $(A, F)$  has a coincidence point. The mappings $A$ and $F$ have a unique  common fixed point if the pair $(A, F)$ is weakly compatible. \\

\noindent \textbf{Corollary 4.2} The conclusions of Theorem 4.1 remain true if we replace the condition (b) with following condition,\\

(b') $\overline{F(X)} \subset B(X) \ or \  \overline{G(X)} \subset A(X).$\\

\noindent \textbf{Theorem 4.2} if  weakly compatible mappings are replaced by any one of the following type mappings (provided remaining settings are same) results of Theorem 4.1 remain true : \\

(1) R-weakly commuting mappings,\\

(2) R-weakly commuting mappings of type $(A_g)$, \\

(3) R-weakly commuting mappings of type $(A_f)$, \\

(4) R-weakly commuting mappings of type $P$, \\

(5) weakly commuting mappings.\\

\noindent \emph{Proof} Provided all the  conditions of Theorem 4.1.1 are satisfied, then there exist  coincidence points for both the pairs. \\

(1) Let $v$ be coincidence point for the pair $(A,F)$, then using R-weak commutativity we get\\

$M(AFv, FAv, t) \geq M(Av, Fv, \frac{t}{R}) = 1$ \\

which amounts to say that $AFV = FAv.$ Thus the pair $(A, F)$ is coincidentally commuting. Similarly $(B, G)$ commutes at all of its coincidence points. Now applying Theorem 4.1, one can see clearly that $A, B, F$ and $G$ have a unique common fixed point. \\

(2) In case $(A, F)$ is an R-weakly commuting pair of type $(A_g),$ then \\

$M(AAv, FAv, t) \geq M(Av, Fv, \frac{t}{R}) = 1 $ \\

which amounts to say that $AAv = FAv.$ Now\\

\begin{eqnarray} \nonumber
M(AFv, FAv, t) &\geq& M( AFv, AAv, \frac{t}{2}) \ast M(AAv, FAv, \frac{t}{2}) \\
&=& 1 \ast 1 = 1 \nonumber
 \end{eqnarray}
this concludes $AFv = FAv.$ \\

Similarly in case of (3), (4) and (5) i.e. if the pair is R-weakly commuting mappings of type $(A_f)$ or type $(P)$ or weakly commuting, then $(A, F)$ also commutes at their  coincidence point. Similarly, we can show that the pair $(B, G)$ is also commuting at coincidence point. Now in settings of Theorem 4.1, all four mappings $ A, B, F$ and $G$ have a unique common fixed point. This completes the proof. \\

\noindent \textbf{Theorem 4.3} Let $\lbrace A_1, A_2, \ldots, A_l \rbrace, \lbrace B_1, B_2, \ldots, B_m \rbrace,\lbrace F_1, F_2 \ldots , F_n \rbrace$ and $\lbrace G_1, G_2, \ldots , G_p \rbrace $ be four finite families of self maps of a fuzzy metric space $(X, M, \ast)$ such that $A = A_1 A_2 \ldots A_l$, $B = B_1 B_2 \ldots B_m$, $F = F_1, F_2 \ldots , F_n$, and $G = G_1 G_2 \ldots G_p $  satisfy all the conditions $(a), (b), (c)$ and $(d)$ of Theorem 4.1 , then \\

\noindent (i) The pair $(A, F)$ has a point of coincidence, \\

\noindent (ii) The pair $(B, G)$ has a point of coincidence, \\

\noindent (iii) If $A_j A_k = A_k A_j, B_r B_s = B_s B_r, F_h F_i = F_i F_h, G_u G_t = G_t G_u, A_j F_h = F_h A_j$ and $B_r G_u = G_u B_r $ for all $j,k \in I_1 = \lbrace 1,2, \ldots , l \rbrace, \ r,s \in I_2 = \lbrace1, 2, \ldots, m\rbrace, \ h,i \in I_3 = \lbrace 1, 2, \ldots, n \rbrace, \ u,t \in I_4 = \lbrace 1,2 \ldots, p\rbrace, then A_j, B_r, F_h \ and \ G_u$ have a common fixed point for all $j \in I_1, r \in I_2, h \in I_3 \ and \ u \in I_4.$ \\

\emph{Proof} Result (i) and (ii) are quick consequence as A, B, F and G satisfy all the conditions of Theorem 4.1. Furthermore mappings A, B, F and G have a unique common fixed point since all the conditions of Theorem 4.1. are satisfied, all we need to prove that z remains the fixed point of all the component mappings. For that we consider \\

\begin{eqnarray} \nonumber
A(A_jz) &=& ((A_1 A_2 \ldots A_l)A_j)z =(A_1 A_2 \ldots A_{l-1})(A_l A_j)z \\ \nonumber
&=& (A_1 \ldots A_{l-1})(A_j A_lz) = (A_1 \ldots A{l-2})(A_{l-1} A_j (A_l z)) \\ \nonumber
&=& (A_1 \ldots A{l-2})( A_j A_{l-1}(A_l z)) = \ldots \\
&=& A_j(A_1 A_2 \ldots A_lz) = A_j(Az) = A_j z. \nonumber
\end{eqnarray}
In the same way, we can prove \\

$A(F_h z) = F_h (Az) = F_h z, \  \ F(F_h z) = F_h(Fz) = F_hz,$\\

$F(A_jz) = A_j(Fz) = A_jz, \ \  B(B_rz) = B_r(Bz) = B_rz,$\\

$B(G_uz) = G_u(Bz) = G_uz,  \ \ G_u(Gz) = G_uz$ \\

and $G(B_rz) = B_r(Gz) = B_rz,$ \\

\noindent one can see clearly that $A_jz$ and $F_hz$ are other fixed points of the the mappings $A$ and $F$ whereas $B_rz$ and $G_uz$ are other fixed points of the mappings $B$ and $G$. For the sake of uniqueness of common fixed points of all four mappings, we get\\

$z = A_jz = B_rz = F_hz = G_uz,$ \\

\noindent which proves that z is a common fixed point of component mappings $A_j, B_r, F_h$ and $G_u.$ \\

\noindent \textbf{Corollary 4.3.} Results of Theorem 4.1 remain true if contractive condition (4.1.1) is interchanged by one of the following contractive conditions. For all $x,y \in X, \delta \in \Psi$ and $\varphi(u$) is an altering distance\\

\noindent (A) $\varphi(M(Fx, Gy, t)) \geq \delta(max\lbrace \varphi(M(Ax, By, t)), \varphi(M(Ax, Fx, t)), \varphi(M(By, Gy, t))\rbrace)$\\
where $\delta : [0,1] \longrightarrow [0,1]$ is an upper semicontinuous function such that $\delta(0) =0$ and $\delta(u) < u$ for $u > 0.$\\

\noindent (B) $\varphi(M(Fx, Gy, t)) \geq k \ (min\lbrace\varphi(M(Ax, By, t)), \varphi(M(Ax, Fx, t)), \varphi(M(By, Gy, t)\rbrace)$ where $k<1$.\\

\noindent (C)  $\varphi(M(Fx, Gy, t)) \geq \delta(\varphi(M(Ax, By, t)), \varphi(M(Ax, Fx, t)), \varphi(M(By, Gy, t))$ 
\noindent where $\delta: [0,1]^3 \to [0,1]$ is an upper semicontinuous function such that max$\lbrace \delta( 0,u,0),\delta(0,0,u),\delta(u,0,0)\rbrace <u  $ for every $u>0.$\\

\noindent (D) $\varphi(M(Fx, Gy, t)) \geq k\varphi(M(Ax, By, t)) - min \lbrace \varphi(M(Ax, Fx, t),  \varphi(M(By, Gy, t)\rbrace$ where $0 < k < 1$. \\

\noindent \emph{Proof} The proof follows from Theorem 4.1 in the settings of Examples 2.1-2.4.\\

\noindent \textbf{Remark 4.1.} Contractive conditions given in Corollary 4.3 with altering distances are new results in the setting of fuzzy metric space. By utilizing the presented altering distances, now we deduce the integral type contractive conditions from our main results, which is also new idea in fuzzy metric space. 

\section{Applications}
\subsection{Application to integral type contractions}

\noindent Branciari proved the following result for integral type contractions as a generalization of Banach fixed point theorem. \\

\textbf{Theorem 5.1.} \cite{Bran} Let (X, d) be a complete metric space and $f : X \longrightarrow X$ be a mapping such that for all $x, y \in X $and $k \in (0, 1)$ \\ 

 $$\int_{0}^{d(fx,fy)}\phi(t)dt \leq k \int_{0}^{d(x,y)} \phi(t)dt > 0,$$  \\

\noindent where $\phi : [0,\infty) \longrightarrow [0,\infty)$ is a Lebesgue measurable mapping (i.e. with finite integral) on each compact subset of $[0,\infty)$ such that for $\epsilon > 0$, $$ \int_{0}^{\epsilon} \phi(t)dt > 0.$$ Then $f$ has a unique fixed point $z \in X$ and for all $x \in X$,
$lim_{n \to \infty} f^nx = z$.\\

Various common fixed point theorems in abstract spaces for compatible, weakly compatible and occasionally weakly compatible mappings satisfying contractions of integral type are proved. Popa and Mocanu \cite{MP} established integral type contractions with the help of altering distances and brought general common fixed point results for integral type inequalities.\\

Next we prove common fixed point theorems for integral type contractive conditions
utilizing new class of altering distance function in the settings of fuzzy metric space which is being introduced by us, which also unifies many results.\\

\noindent \textbf{Lemma 5.1.} The function $ \varphi(s) = \int_{0}^{1-s} \phi(x)dx,$ where $\varphi(s)$ is an altering distance function and $\phi(x)$ is a Lebesgue measurable function  .  \\

We next prove a common fixed point theorem for pairs of self mappings satisfying integral type contraction .\\

\noindent \textbf{Theorem 5.2.} Let $A, B, F \ and \ G$ be self mappings of a fuzzy metric space $(X, M, \ast)$. Suppose that \\

\noindent (a) Either the pair (A, F) or (B, G) share the property (E.A.),\\
(b) Either B(X) $\subset$ F(X) or G(X) $\subset$ A(X),\\
(c) A(X) is a closed subset of X(or B(X) is  closed subset of X) and \\

 $$ \psi  \Big( \int_{0}^{1-M(Fx, Gy, t)} \phi(s)ds, \int_{0}^{1-M(Ax, By, t)}\phi(s)ds, \int_{0}^{1-M(Ax, Fx, t)}\phi(s)ds, \int_{0}^{1-M(By, Gy, t)}\phi(s)ds \Big) \geq 0.$$ (5.1.1) \\
 
 for all $x,y \in X, \ \psi \in \Psi, \ \phi \in \Phi $(as mentioned in Theorem 5.1). If the self mapping pairs are weakly compatible, then there exists a unique common fixed point for mappings $A, B, F$ and $G$.\\

\noindent \emph{Proof} From Lemma 5.1, we have \\

\noindent  $\varphi(M(Fx, Gy, t)) = \int_{0}^{1- M(Fx, Gy, t)} \phi(s)ds, \  \varphi(M(Ax, By, t)) = \int_{0}^{1- M(Ax, By, t)}\phi(s)ds, \  \varphi(M(Ax, Fx, t)) =\int_{0}^{1-M(Ax, Fx, t)}\phi(s)ds \ and \ \varphi(M(By, Gy, t))=\int_{0}^{1-M(By, Gy, t)}\phi(s)ds$. \\

\noindent Then by contractive condition (4.1.1), we have \\

$\psi \Big( \varphi(M(Fx, Gy, t)), \varphi(M(Ax, By, t)), \varphi(M(Ax, Fx, t)), \varphi(M(By, Gy, t)) \Big) \geq 0.$\\

\noindent We can see $\varphi(s)$ is an altering distance function from Lemma 5.1. Results of Theorem 5.1 follow from Theorem 4.1 since all the conditions of Theorem 4.1 are satisfied. Hence the theorem is proved.\\

\noindent \textbf{Corollary 5.1.} The results of Theorem 5.1 will not change if contractive condition (5.1.1) is replaced by one of the following inequalities (for all $x,y \in X, \psi \in \Psi$ and $\phi \in \Phi$) \\

\begin{eqnarray} \nonumber
(A) \int_0^{1-M(Fx, Gy, t)} \phi(u)du \geq a \ max \Big\lbrace \int_0^{1-M(Ax, By, t)} \phi(u)du, 
 \int_0^{1-M(Ax, Fx, t)} \phi(u)du,\\ \int_0^{1-M(By, Gy, t)}\phi(u)du \Big\rbrace \nonumber
   \end{eqnarray} 
where $0 \leq a < 1$ and $\phi \in \Phi$. \\

\begin{eqnarray} \nonumber
(B) \int_0^{1- M(Fx, Gy, t)} \phi(u) du \geq \delta \Big(max\Big\lbrace \int_0^{1-M(Ax, By, t)}\phi(u)du ,\int_0^{1-M(Ax, Fx, t)}\phi(u)du, \\ \int_0^{1-M(By, Gy, t)}\phi(u)du\Big\rbrace \Big). \nonumber
\end{eqnarray}
\noindent where $\delta:  \mathfrak{R}_+ \longrightarrow \mathfrak{R}_+ $ is continuous function such that $\delta(u) < u$ for all $u > 0$ and $\phi \in \Phi$.\\

\noindent \textbf{Proof.} We can see the proof follows from Theorem 5.2 in the settings of Examples 2.5-2.6.\\

\noindent \textbf{Remark 5.1.}  Theorem 2 and Theorem 3 of  Murthy et al \cite{PPM} can be deduced from the corollary 5.1. Hence our results generalize and unify existing results. 

\subsection{Application to system of functional equation}
Let S and T be Banach spaces, W $\subset$ S be a state space and D $\subset$ T be a decision space. Now as an application of Theorem 4.1 we discuss the solvability of following system of functional equations arising in dynamic programming: \\
\begin{equation}
\begin{cases}
P_i(x)  =  \sup_{y\in D} \lbrace q(x,y) + L_i(x,y, P_i(\tau (x,y))) \rbrace, x \in W, i = 1,2,\\ 
Q_i(x)  =  \sup_{y\in D} \lbrace q(x,y) + N_i(x,y, P_i(\tau (x,y))) \rbrace, x \in W, i = 1,2,
\end{cases}
\end{equation}
where $\tau : W \times D \longrightarrow W, q : W\times D \longrightarrow \mathbb{R}, L_i, N_i : W \times D \times \mathbb{R}\longrightarrow \mathbb{R}. $ \\ 
Let C(W) represents the space of all real valued continuous functions on W. Precisely, this space will have the metric given by\\

$d(r,p) = \sup_{x \in W} \vert r(x) - p(x) \vert,$ for all $r,p \in C(W)$ \\

is a complete metric space. \\
We will prove the following result.\\

\textbf{Theorem 5.3}
Let $L_i, N_i: W \times D \times \mathbb{R},$ for $i = 1,2,$ be bounded functions and let $U_i, V_i : C(W) \longrightarrow C(W), for i =1,2,$ be four operators defined as, \\
\begin{eqnarray}
U_ir(x)  =  \sup_{y\in D} \lbrace q(x,y) + L_i(x,y, r(\tau (x,y))) \rbrace, x \in W, i = 1,2,\\ 
V_ip(x)  =  \sup_{y\in D} \lbrace q(x,y) + N_i(x,y, p(\tau (x,y))) \rbrace, x \in W, i = 1,2, \\ \nonumber
\end{eqnarray} 
for all $r,p \in C(W)$ and $x \in W.$ Assuming that the following conditions hold:\\

(i) there exists $\{r_n\} \in C(W)$ such that $\lim_{n\to \infty} U_1r_n = \lim_{n\to \infty}U_2r_n = r^\ast \in C(W)$ \\ and $\lim_{n\to \infty} \sup_{x\in W} \vert U_1U_2r_n - U_2U_1r_n \vert =0, $

(ii) there exists $\{p_n\} \in C(W)$ such that $\lim_{n\to \infty} V_1p_n = \lim_{n\to \infty} V_2p_n = p^\ast \in C(W)$ \\ and $\lim_{n\to \infty} \sup_{x\in W} \vert V_1V_2p_n - V_2V_1p_n \vert =0, $ \\

(iii) $\vert L_1(x,y,r(\tau(x,y))) - N_1 (x,y,p(\tau(x,y)))\vert \leq \Theta (r,p)$\\

where $\Theta(r,p)= \lambda( max\{ \varphi(d(U_2r, V_2p)), \varphi(d(U_2r, U_1r)), \varphi(d(V_2p,V_1p)\})$ with $\lambda(x) \geq x$ \\
 with $\varphi(t) = t-1$. Then the system of functional equations (1) has unique bounded solution. \\
 
\noindent \emph{Proof}  The system of functional equation (1) will have a unique bounded solution if and only if the operators in (2) and (3) have unique common fixed point. Now since $L_1, L_2, N_1$ and $N_2$ are bounded, there exists a positive number $\Lambda  $ such that, \\
\begin{eqnarray}
 sup\{\vert L_i(x,y,z)\vert, \vert N_i(x,y,z)\vert : (x,y,z) \in W \times D \times \mathbb{R}, i= 1,2\} \leq \Lambda. \nonumber
\end{eqnarray}
Let $\varepsilon$ be an arbitrary positive number, $x \in W$ and $r_1, r_2 \in C(W),$ then there exist $y_1,y_2 \in D $ such that \\

\begin{eqnarray}
U_1r_1(x) &<& q(x,y_1) + L_1(x,y_1,r_1(\tau(x,y_1))) + \varepsilon,\\
V_1r_2(x) &<& q(x,y_2) + N_1(x,y_2,r_2(\tau(x,y_2))) + \varepsilon,\\
U_1r_1(x) &\geq & q(x,y_2) + L_1(x,y_2,r_1(\tau(x,y_2))),\\
V_1r_2(x) &\geq & q(x,y_1) + N_1(x, y_1,r_2(\tau(x,y_1))).
\end{eqnarray}
Now from (4) and (7) we get, \\
\begin{eqnarray}
U_1r_1(x) - V_1r_2(x) &<& L_1(x,y_1,r_1(\tau(x,y_1))) - N_1(x, y_1,r_2(\tau(x,y_1))) + \varepsilon \nonumber \\
&<& \vert L_1(x,y_1,r_1(\tau(x,y_1))) - N_1(x, y_1,r_2(\tau(x,y_1)))\vert  + \varepsilon \\
\text{Similarly from (5) and (6)} \nonumber \\
V_1r_2(x) - U_1r_1(x) &<&  N_1(x,y_2,r_2(\tau(x,y_2))) -  L_1(x,y_2,r_1(\tau(x,y_2))) + \varepsilon \nonumber \\
&\leq & \vert N_1(x,y_2,r_2(\tau(x,y_2))) -  L_1(x,y_2,r_1(\tau(x,y_2)))\vert + \varepsilon \\
\text{From (8) and (9), we get,} \nonumber \\ 
\vert U_1r_1(x)-V_1r_2(x)\vert &\leq & \vert L_1(x,y_1,r_1(\tau(x,y_1)))-N_1(x,y_1,r_2(\tau(x,y_1)))\vert  +\varepsilon \\
d(U_1r_1,V_1r_2) &\leq & \vert L_1(x,y_1,r_1(\tau(x,y_1)))-N_1(x,y_1,r_2(\tau(x,y_1)))\vert  +\varepsilon
\end{eqnarray}
Since the above inequality does not depend on $x \in W$  and positive number $\varepsilon$ is taken arbitrary, then we can conclude \\

$$d(U_1r_1,V_1r_2) \leq  \vert L_1(x,y_1,r_1(\tau(x,y_1)))-N_1(x,y_1,r_2(\tau(x,y_1)))\vert$$ \\
$$\varphi(d(U_1r_1,V_1r_2)) \leq  \vert L_1(x,y_1,r_1(\tau(x,y_1)))-N_1(x,y_1,r_2(\tau(x,y_1)))\vert$$ \\
where $\varphi$ is an altering distance function such that $\varphi(t) = t-1$\\

now from assumption (iii),

$$\varphi(d(U_1r_1,V_1r_2)) \leq  \Theta (r,p)$$
where $\Theta(r_1,r_2)= \lambda( max\{ \varphi(d(U_2r_1, V_2r_2)), \varphi(d(U_2r_1, U_1r_1)), \varphi(d(V_2r_2,V_1r_2)\})$ \\

and so Theorem 4.1 is applicable with $\psi(u_1,u_2,u_3,u_4) = \lambda(max\lbrace u_2,u_3,u_4 \rbrace) - u_1 $ where $\lambda(u) > u$. Consequently the mapping T has a unique fixed point, that is, the system of functional equations (1) has a unique bounded solution.
\section{Illustrative Example} Let $X = [0,1]$ be the set of all real numbers with usual metric d defined by $d(x,y) = \mid x-y \mid$ for all $x,y \in X$. Define \\

$M(x,y,t) = \dfrac{t}{t+d(x,y)}$, for all $x,y \in X, t > 0. $ \\

\noindent Then $(X, M, \ast)$ is a fuzzy metric space. Define $Ax = \frac{x}{2}, Fx = x, Bx = \frac{x}{4}, Gx = 0 \ \forall x \in X.$  The pair $(A, F)$ is weakly compatible since $AF(0) = FA(0) = \lbrace 0 \rbrace.$ Also $ A(X) = [0, \frac{1}{2}]$ is closed in X and the pair $(A, F)$ satisfies property(E.A.) $(consider \lbrace x_n\rbrace = \frac{1}{n})$ as \\

$\lim_{n \to \infty} Ax_n = \lim_{n \to \infty} Fx_n =\lbrace 0 \rbrace.$ \\

\noindent First we verify the condition(4.1.1) of Theorem 4.1. We set the mappings $\psi$  as \\

$\psi (u_1, u_2, u_3, u_4) = u_1 - k min\lbrace u_2, u_3, u_4 )$, where $ 0<k<1$\\

\noindent we have to prove the corresponding contractive condition \\

\noindent $\varphi(M(Fx, Gy, t)) \geq k  min\lbrace \varphi(M(Ax, By, t)), \varphi(M(Ax, Fx, t), \varphi(M(By, Gy, t)).$\\ 

where $0<k<1$\\

\noindent Suppose the altering distance $\varphi (\lambda) = 1 - \lambda $  where  $ \lambda \in [0, 1]$, \\

 If $x,y \in [0, 1]$ and $ y > 4x$ \\

\begin{eqnarray}\nonumber
min \Big\lbrace \varphi(M(Ax, By, t)), \varphi(M(Ax, Fx, t), \varphi(M(By, Gy, t))\Big \rbrace &=& min \Big\lbrace  \dfrac{\frac{y}{4}- \frac{x}{2}}{t + \frac{y}{4}- \frac{x}{2} }, \dfrac{\frac{x}{2}}{t + \frac{x}{2}}, \dfrac{\frac{y}{4}}{t + \frac{y}{4}}\Big\rbrace \\
&=& \dfrac{\frac{x}{2}}{t + \frac{x}{2}} \nonumber 
\end{eqnarray}
Now, \begin{eqnarray} \nonumber
\varphi(M(Fx, Gy, t)) &=& \dfrac{x}{t + x} \\ \nonumber
&\geq& k \ \dfrac{x}{t + x} (where\ 0<k<1)\\ \nonumber 
&\geq& k \ \dfrac{\frac{x}{2}}{t + \frac{x}{2}}\\ \nonumber
&\geq& k \ min \Big\lbrace \varphi(M(Ax, By, t)), \varphi(M(Ax, Fx, t), \varphi(M(By, Gy, t))\Big \rbrace .
\end{eqnarray}
Therefore all the conditions of Theorem 4.1 are satisfied. Here, 0 is the coincidence as well as unique common fixed point of the mappings A, B, F and G. This example verifies our main result.

\end{document}